\documentstyle[12pt]{article}
\begin{document}
\title {Approximate $(\sigma-\tau)$-Contractibility \footnote{{\it 2000 Mathematics Subject Classification}. 39B82; 46H25; 39B52; 47B47.\\
{\it Key words and phrases}. Hyers--Ulam stability, $(\sigma-\tau)$-derivation, inner $(\sigma-\tau)$-derivation, approximate $(\sigma-\tau)$-amenability, approximate $(\sigma-\tau)$-contractibility.}}
\author{{\bf Mohammad Sal Moslehian} \\ Dept. of Math., Ferdowsi Univ.\\ P. O. Box 1159, Mashhad 91775\\ Iran\\ E-mail: msalm@math.um.ac.ir\\Home: http://www.um.ac.ir/$\sim$moslehian/}
\date{}
\maketitle
\begin{abstract}
In this paper, the generalized Hyers--Ulam--Rassias stability of $(\sigma-\tau)$-derivations on normed algebras into Banach bimodules is established. We introduce the notion of approximate $(\sigma-\tau)$-contractibility and prove that a Banach algebra is $(\sigma-\tau)$-contractible if and only if it is approximately $(\sigma-\tau)$-Contractible.
\end{abstract}
\newpage

\section{Introduction.}

The stability theory of functional equations was started in 1940 with a problem on approximate homomorphisms raised by S. M. Ulam; cf. $\cite{ULA}$. Indeed, as noted by Z. P\' ales $\cite{PAL}$, the first stability theorem was discovered in 1925 anticipating the question of Ulam $\cite{P-S}$. In 1941, D. H. Hyers gave a partial solution of Ulam's problem for linear mappings. Since then many mathematicians have been working on this area of research; cf. $\cite{CZE}$ and $\cite{FOR}$. The stability of derivations was studied by C.-G. Park $\cite{PAR1}, \cite{PAR2}$ and by the author $\cite{MOS}$.

Let $A$ be an algebra and $X$ be an $A$-bimodule. A linear mapping $d:A\to X$ is called a $(\sigma-\tau)$-derivation if there exists two linear operators $\sigma$ and $\tau$ on $A$ such that $d(ab)=d(a)\sigma(b)+\tau(a)d(b)$ for all $a,b\in A$. Familiar examples are\\
(i) the ordinary derivations from $A$ to $X$;\\
(ii) the so-called $(\sigma-\tau)$-inner derivations i.e. those are defined by $d_x(a)=x\sigma(a)-\tau(a)x$ for a fixed arbitrary element $x\in X$ and endomorphisms $\sigma$ and $\tau$ on $A$. Note that $d_x(ab)=d_x(a)\sigma(b)+\tau(a)d_x(b)$.\\
(iii) all endomorphisms $\phi$ on $A$. Notice that such $\phi$ is clearly a $(\frac{1}{2}\phi-\frac{1}{2}\phi)$-derivation. See $\cite{A-N}$.

The present paper is devoted to study of stability of $(\sigma-\tau)$-derivations and to generalize some results in $\cite{PAR1}$. We also introduce $(\sigma-\tau)$-contractibility ($(\sigma-\tau)$-amenability) and approximate $(\sigma-\tau)$-contractibility (approximate $(\sigma-\tau)$-amenability) and investigate their relationships which may seem to be interesting on their own right.

\section{Stability of $(\sigma-\tau)$-Derivations.}

Throughout the section, $A$ denotes a (not necessary unital) normed algebra and $X$ is a Banach $A$-bimodule. Our aim is to establish the generalized Hyers--Ulam--Rassias stability of $(\sigma-\tau)$-derivations. We indeed extend main results of C.-G. Park $\cite{PAR1}$ to $(\sigma-\tau)$-derivations. We use the direct method which was first devised by D. H. Hyers to construct an additive function from an approximate one; cf. $\cite{HYE}$ and G\u avruta's technique $\cite{GAV}$. Some ideas of $\cite{PAR2}$ are also applied.

{\bf 2.1. Theorem.} Suppose $f:A\to X$ is a mapping with $f(0)=0$ for which there exist maps $g_1,g_2:A\to A$ with $g_1(0)=g_2(0)=0$ and a function $\varphi: A\times A\to [0,\infty)$ such that 
\begin{eqnarray}
\widetilde{\varphi}(a,b):=\frac{1}{2}\displaystyle{\sum_{n=0}^\infty}2^{-n}\varphi(2^na,2^nb)<\infty
\end{eqnarray}
\begin{eqnarray}
\|f(\lambda a+\lambda b)-\lambda f(a)-\lambda f(b)\|\leq \varphi(a,b)
\end{eqnarray}
\begin{eqnarray}
\|g_k(\lambda a+\lambda b)-\lambda g_k(a)-\lambda g_k(b)\|\leq \varphi(a,b)
\end{eqnarray}
\begin{eqnarray}
\|f(ab)-f(a)g_1(b)-g_2(a)f(b)\|\leq \varphi(a,b)
\end{eqnarray}
for $k=1,2$, for all $\lambda\in T=\{\lambda\in {\bf C}: |\lambda|=1\}$ and for all $a,b\in A$. Then there exist unique linear maps $\sigma$ and $\tau$ from $A$ to $X$ satisfying $\|g_1(a)-\sigma(a)\|\leq \widetilde{\varphi}(a,a)$ and $\|g_2(a)-\tau(a)\|\leq \widetilde{\varphi}(a,a)$, and there exists a unique $(\sigma-\tau)$-derivation $d:A\to X$ such that 
\begin{eqnarray}
\|f(a)-d(a)\|\leq\widetilde{\varphi}(a,a)
\end{eqnarray}
for all $a\in A$. Moreover,\\
(i) if
\begin{eqnarray}
\|g_2(ab)-g_2(a)g_2(b)\|\leq \varphi(a,b)
\end{eqnarray}
for all $a,b\in A$, and either $A$ has no nonzero divisor of zero or $d$ is surjective and ${\rm ran}(A)=\{0\}$, then either $d=0$ or both linear mappings $\sigma$ and $\tau$ are endomorphisms on $A$.\\
(ii) if $f$ and $\psi(a)=\widetilde{\varphi}(a,a)$ are continuous at a point $a_0$ then $d$ is continuous on $A$.

{\bf Proof.} Putting $\lambda=1$ in $(2)$, we have 
\begin{eqnarray}
\|f(a+b)-f(a)-f(b)\|\leq\varphi(a,b)~~~ a,b\in A
\end{eqnarray}
Now we use the G\u avruta method on inequality $(4)$ (see $\cite{GAV}$ and $\cite{RAS}$). One can use the induction to show that
\begin{eqnarray}
\|\frac{f(2^na)}{2^n}-f(a)\|\leq\frac{1}{2}\displaystyle{\sum_{k=0}^{n-1}}2^{-k}\varphi(2^ka,2^ka)
\end{eqnarray}
and $$\|\frac{f(2^pa)}{2^p}-\frac{f(2^qa)}{2^q}\|\leq\frac{1}{2}\displaystyle{\sum_{k=q}^{p-1}}2^{-k}\varphi(2^ka,2^ka)$$
for all $a\in A$, all $n$ and all $p>q$. It follows from the convergence $(1)$ that the sequence $\{\frac{f(2^na)}{2^n}\}$ is Cauchy. Due to the completeness of $X$, this sequence is convergent. Set 
\begin{eqnarray}
d(a):= \displaystyle{\lim_{n\to\infty}}\frac{f(2^na)}{2^n}
\end{eqnarray}
Applying $(2)$, we get $\|2^{-n}f(2^n(\lambda a+\lambda b))-2^{-n}\lambda f(2^na)-2^{-n}\lambda f(b)\|\leq 2^{-n}\varphi(2^na,2^na)$. Passing to the limit as $n\to\infty$ we obtain $d(\lambda a+\lambda b)=\lambda d(a)+\lambda d(b)$. Note that the convergence of $(1)$ implies that $$\displaystyle{\lim_{n\to\infty}}2^{-n}\varphi(2^na,2^na)=0.$$
Next, let $\gamma\in {\bf C} (\gamma\neq 0)$ and let $M$ be a positive integer greater than $4|\gamma|$.
Then $|\frac{\gamma}{M}|<\frac{1}{4}<1-\frac{2}{3}=1/3$. By Theorem 1 of $\cite{K-P}$, there exist three elements $\theta_1,
\theta_2, \theta_3\in T$ such that $3\frac{\gamma}{M}=\theta_1+\theta_2+\theta_3$. By the additivity of $d$
we get $d(\frac{1}{3}a)=\frac{1}{3}d(a)$ for all $ a\in A$. Therefore,
\begin{eqnarray*}
d(\gamma a)& = & d(\frac{M}{3}\cdot 3 \cdot \frac{\gamma}{M}a)=Md(\frac{1}{3}\cdot 3\cdot \frac{\gamma}{M}a)
=\frac{M}{3}d(3\cdot \frac{\gamma}{M}a)\\ & = & \frac{M}{3}d(\theta_1a+\theta_2a+\theta_3a)
=\frac{M}{3}(d(\theta_1a)+d(\theta_2a)+d(\theta_3a)) \\ & = & \frac{M}{3}(\theta_1+\theta_2+\theta_3)d(a)
=\frac{M}{3}\cdot 3\cdot \frac{\gamma}{M}=\gamma d(a)
\end{eqnarray*}
for all $a\in A$.\\
Thus $d$ is (${\bf C}$-)linear.

Moreover, it follows from $(8)$ and $(9)$ that $\|f(a)-d(a)\|\leq\widetilde{\varphi}(a,a)$ for all $a\in A$. It is known that additive mapping $d$ satisfying $(2)$ is unique $\cite{GAV}$.

Similarly one can use $(3)$ to show that there exist unique linear mappings $\sigma$ and $\tau$ defined by $\displaystyle{\lim_{n\to\infty}}2^{-n}g_1(2^nb)$ and $\displaystyle{\lim_{n\to\infty}}2^{-n}g_2(2^na)$, respectively.

Replacing $a, b$ in $(4)$ by $2^na$ and $2^nb$ respectively, we have
\begin{eqnarray*}
\|f(2^{2n}ab)-f(2^na)g_1(2^nb)-g_2(2^na)f(2^nb)\|\leq\varphi(2^na,2^nb)
\end{eqnarray*}
\begin{eqnarray*}
\|2^{-2n}f(2^{2n}ab)-2^{-n}f(2^na)2^{-n}g_1(2^nb)-2^{-n}g_2(2^na)2^{-n}f(2^nb)\|\leq2^{-2n}\varphi(2^na,2^nb)
\end{eqnarray*}
It follows that
\begin{eqnarray*}
d(ab)=d(a)\sigma(b)+\tau(a)d(b).
\end{eqnarray*}
Therefore $d$ is a $(\sigma-\tau)$-derivation.

(i) Inequality $(6)$ yields
\begin{eqnarray*}
\|2^{-2n}g_2(2^{2n}ab)-2^{-n}g_2(2^na)2^{-n}g_2(2^nb)\|\leq 2^{-n}\varphi(2^na,2^nb).
\end{eqnarray*}
for all $a,b\in A$. So that $\tau(ab)=\tau(a)\tau(b)$ for all $a,b\in A$. \\
On the other hand 
\begin{eqnarray*}
d(cab)&=&d(c)\sigma(ab)+\tau(c)d(ab)\\
d(c)\sigma(ab)&=&d(cab)-\tau(c)d(ab)\\
&=&(d(ca)\sigma(b)+\tau(ca)d(b))-\tau(c)d(ab)\\
&=&(d(c)\sigma(a)+\tau(c)d(a))\sigma(b)+\tau(ca)d(b)-\tau(c)d(ab)\\
&=&d(c)\sigma(a)\sigma(b)+\tau(c)d(a)\sigma(b)+\tau(ca)d(b)-\tau(c)d(ab)\\
&=&d(c)\sigma(a)\sigma(b)+\tau(c)d(a)\sigma(b)+\tau(ca)d(b)\\
&&-\tau(c)(d(a)\sigma(b)+\tau(a)d(b))\\
d(c)(\sigma(ab)-\sigma(a)\sigma(b))&=&(\tau(ca)-\tau(c)\tau(a))d(b)=0
\end{eqnarray*}
Hence $d(c)(\sigma(ab)-\sigma(a)\sigma(b))=0$. Therefore if $d\neq 0$, under any one of the assumptions that $A$ has no nonzero divisor of zero, or $d$ is surjective and ${\rm ran}(A)=\{0\}$ we conclude that $\sigma$ is an endomorphism.

(ii) If $d$ were not continuous at a point $a\in A$ then there would be an integer $m$ and a sequence $\{a_n\}$ of $A$ such that $\displaystyle{\lim_{n\to\infty}}a_n=0$ and $\|d(a_n)\|>\frac{1}{m}$ for all $n$. Let $k$ be an integer greater than $m(2\psi(a_0)+1)$. Since $\displaystyle{\lim_{n\to\infty}}f(ka_n+a_0)=f(a_0)$, there is an integer $N$ such that for all $n\geq N, \|f(ka_n+a_0)=f(a_0)\|<1$. Hence 
\begin{eqnarray*}
2\psi(a_0)+1&<&\frac{k}{m}\\
&<&\|d(ka_n)\|=\|d(ka_n+a_0)-d(a_0)\|\\
&\leq&\|d(ka_n-a_0)-f(ka_n+a_0)\|+\|f(ka_n+a_0)-f(a_0)\|\\
&&+\|f(a_0)-d(a_0)\|\\
&\leq&\psi(ka_n+a_0)+1+\psi(a_0)
\end{eqnarray*}
for all $n\geq N$. Tending $n$ to $\infty$ we conclude that $2\psi(a_0)+1<\frac{k}{m}\leq 2\psi(a_0)+1$, a contradiction (see $\cite{HYE}$).$\Box$

{\bf 2.2. Remark.} 
(i) The algebra ${\bf C}[x]$ consisting of all complex polynomials in one variable has no zero-divizor.\\
(ii) ${\rm ran}A={\rm lan}A=0$ for each normed algebra $A$ with a bounded approximate identity or a unital normed algebra $A$ containing an invertible element.\\
(iii) A similar statement to Theorem 2.1 holds if 
\begin{eqnarray*}
\|g_1(ab)-g_1(a)g_1(b)\|\leq \varphi(a,b)
\end{eqnarray*}
for all $a,b\in A$, and either $A$ has no nonzero divisor of zero or $d$ is surjective and ${\rm lan}(A)=\{0\}$.

{\bf 2.3. Proposition.} Suppose the normed algebra $A$ is spanned by a set $S\subseteq A$, $f:A\to X$ is a mapping with $f(0)=0$ for which there exist maps $g_1,g_2:A\to A$ with $g_1(0)=g_2(0)=0$ and a function $\varphi: A\times A\to [0,\infty)$ such that 
\begin{eqnarray*}
\widetilde{\varphi}(a,b):=\frac{1}{2}\displaystyle{\sum_{n=0}^\infty}2^{-n}\varphi(2^na,2^nb)<\infty
\end{eqnarray*}
\begin{eqnarray}
\|f(\lambda a+\lambda b)-\lambda f(a)-\lambda f(b)\|\leq \varphi(a,b)
\end{eqnarray}
\begin{eqnarray*}
\|g_k(\lambda a+\lambda b)-\lambda g_k(a)-\lambda g_k(b)\|\leq \varphi(a,b)
\end{eqnarray*}
for $k=1,2$, for $\lambda=1,{\rm i}$ and for all $a,b\in A$; and
\begin{eqnarray*}
\|f(ab)-f(a)g_1(b)-g_2(a)f(b)\|\leq \varphi(a,b)
\end{eqnarray*}
for all $a,b\in A$. If for $k=1,2$ and each fixed $a\in A$ the functions $t\mapsto f(ta)$, $t\mapsto g_k(ta)$ are continuous on ${\bf R}$ then there exist unique linear maps $\sigma$ and $\tau$ from $A$ to $X$ satisfying $\|g_1(a)-\sigma(a)\|\leq \widetilde{\varphi}(a,a)$ and $\|g_2(a)-\tau(a)\|\leq \widetilde{\varphi}(a,a)$, and a unique $(\sigma-\tau)$-derivation $d:A\to X$ such that 
\begin{eqnarray*}
\|f(a)-d(a)\|\leq\widetilde{\varphi}(a,a)
\end{eqnarray*}
for all $a\in A$.

{\bf Proof.} Put $\lambda=1$ in (10). It follows from the proof of Theorem 2.1 that there exists a unique additive mapping $d:A\to X$ given by $d(a)=\displaystyle{\lim_{n\to\infty}}\frac{f(2^na)}{2^n}, a\in A$. By the same reasoning as in the proof of the theorem of $\cite{RAS}$, the mapping $d$ is ${\bf R}$-linear. 

Assuming $b=0$ and $\lambda={\rm i}$ in $(10)$ we have $\|f({\rm i}a)-{\rm i}f(a)\|\leq\varphi(a,0), a\in A$. Hence $\frac{1}{2^n}\|f(2^n{\rm i}a)-{\rm i}f(2^na)\|\leq\varphi(2^na,0)a\in A$.  The right hand side tends to zero as $n\to \infty$, so that $d({\rm i}a)=\displaystyle{\lim_{n\to\infty}}\frac{f(2^n{\rm i}a)}{2^n}=\displaystyle{\lim_{n\to\infty}}\frac{if(2^na)}{2^n}={\rm i}d(a), a\in A$.\\
For every $\lambda\in {\bf C}, \lambda=r_1+{\rm i}r_2$ in which $r_1,r_2\in{\bf R}$. Hence $d(\lambda a)=d(r_1a+{\rm i}r_2a)=r_1d(a)+{\rm i}r_2d(a)=\lambda d(a)$. Thus $d$ is ${\bf C}$-linear.

Similarly there exist linear maps $\sigma$ and $\tau$ from $A$ to $X$ satisfying the required inequalities. One can easily verify $\|f(a)-d(a)\|\leq\widetilde{\varphi}(a,a)$ and $d(ab)=d(a)\sigma(b)+\tau(a)d(b)$. Since $A$ is linearly generated by $S$ we conclude that $d$ is $(\sigma-\tau)$-derivation.$\Box$

{\bf 2.4. Remark.} (i) Similar statements to Theorem 2.1 and Proposition 2.3 hold if $\varphi(a,b)=\alpha+\beta(\|a\|^p+\|b\|^p)$. Note that
\begin{eqnarray*}
\widetilde{\varphi}(a,b):=\frac{1}{2}\displaystyle{\sum_{n=0}^\infty}2^{-n}\varphi(2^na,2^nb)=\alpha+\beta\frac{\|a\|^p+\|b\|^p}{2-2^p}.
\end{eqnarray*}
(ii) Conclusion of Proposition 2.3 holds if $A$ is a $C^*$-algebra and $S$ is the unitary group of $A$ or the positive part of $A$.

\section{Approximately $(\sigma-\tau)$-Contractible}

Throughout this section $A$ is a Banach algebra and $X$ is a Banach $A$-bimodule. In addition, $\sigma$ and $\tau$ are assumed to be bounded endomorphisms on $A$.

If every bounded $(\sigma-\tau)$-derivation is $(\sigma-\tau)$-inner, then $A$ is said to be $(\sigma-\tau)$-contractible. In particular, the ordinary contractibility is indeed $({\rm id}-{\rm id})$-contractibility where ${\rm id}$ denotes the identity map.

We say a mapping $f:A\to X$ with $f(0)$=0 is approximate $(\sigma-\tau)$-derivation if it is continuous at a point $a\in A$ and there is a positive number $\alpha>0$ such that $\|f(\lambda a+\lambda b)-\lambda f(a)-\lambda f(b)\|\leq \alpha$ and $\|f(ab)-f(a)\sigma(b)-\tau(a)f(b)\|\leq\alpha$ for all $\lambda\in T$ and for all $a,b\in A$.

The Banach algebra $A$ is called approximately $(\sigma-\tau)$-contractible if for every approximate $(\sigma-\tau)$-derivation there exist a positive number $\beta$ and an $x\in X$ such that $\|x\sigma(a)-\tau(a)x-f(a)\|\leq\beta$.

{\bf 3.1. Theorem.} A Banach algebra $A$ is approximately $(\sigma-\tau)$-contractible if and only if $A$ is $(\sigma-\tau)$-contractible.

{\bf Proof.} Let $A$ be $(\sigma-\tau)$-contractible, $\alpha>0$ and $f:A\to X$ is a mapping  which is continuous at a point $a\in A$ and $f(0)=0$ such that
\begin{eqnarray*}
\|f(\lambda a+\lambda b)-\lambda f(a)-\lambda f(b)\|\leq \alpha
\end{eqnarray*}
\begin{eqnarray*}
\|f(ab)-f(a)\sigma(b)-\tau(a)f(b)\|\leq \alpha
\end{eqnarray*}
for all $\lambda\in T$ and for all $a,b\in A$.\\
By Theorem 2.1 there exists a bounded $(\sigma-\tau)$-derivation $d:A\to X$ defined by $d(a)=\displaystyle{\lim_{n\to\infty}}\frac{f(2^na)}{2^n}, a\in A$ which satisfies $\|d(a)-f(a)\|\leq \alpha$. Since $A$ is $(\sigma-\tau)$-contractible, there is some $x\in X$ such that $d(a)=x\sigma(a)-\tau(a)x$. Hence 
\begin{eqnarray*}
\|x\sigma(a)-\tau(a)x-f(a)\|\leq \|x\sigma(a)-\tau(a)x-d(a)\|+\|d(a)-f(a)\|\leq \alpha
\end{eqnarray*}
Therefore $A$ is approximately $(\sigma-\tau)$-contractible.

Conversely, let $A$ be approximately $(\sigma-\tau)$-contractible and $d:A\to X$ be a bounded $(\sigma-\tau)$-derivation. Then $d$ is trivially an approximate $(\sigma-\tau)$-derivation. Due to the approximate $(\sigma-\tau)$-contractibility of $A$, there exist $\beta>0$ and $x\in X$ such that $\|x\sigma(a)-\tau(a)x-d(a)\|\leq \beta$. Replacing $a$ by $2^na$ in the later inequality we get
$\|x\sigma(a)-\tau(a)x-d(a)\|\leq 2^{-n}\beta$ for all positive integer $n$. Hence $x\sigma(a)-\tau(a)x=d(a)$. It follows that $A$ is $(\sigma-\tau)$-contractible.$\Box$

One can similarly define the notions of $(\sigma-\tau)$-amenability and approximate $(\sigma-\tau)$-amenability and establish the following theorem.

{\bf 3.2. Theorem.} A Banach algebra $A$ is approximately $(\sigma-\tau)$-amenable if and only if $A$ is $(\sigma-\tau)$-amenable.

\end{document}